\documentclass[a4paper,12pt]{article}

\usepackage{setspace}
\usepackage[active]{srcltx}
\usepackage{makeidx}
\usepackage{multicol}
\usepackage{latexsym}
\usepackage{amsmath}
\usepackage{amssymb}
\usepackage{amscd}
\usepackage{amsthm}
\usepackage[english]{babel}
\usepackage{xypic}

\newtheorem{theorem}{Theorem}[section] 
\newtheorem{proposition}[theorem]{Proposition}
\newtheorem{lemma}[theorem]{Lemma}

\newtheorem{definition}{Definition}[section] 
\newtheorem{example}{Example}[section]

\numberwithin{equation}{subsection}
\newtheorem*{Teo}{Theorem \ref{PeetreThm}}
\newtheorem*{TeoNumerado}{Theorem \ref{NonLinearPeetre}}

\def\proof{\noindent \textit{Proof: }}

\def\qed{\hfill $\square$}

\def\O{\mathcal{O}}

\def\proof{\medskip \noindent \textit{Proof:} }

\def\reg{\mathrm{Hom} \hskip -.2cm \phantom{a}_{reg}}

\def\Sdiff{\underset{^{\sim \! \! \sim \! \! \sim}}{\mathrm{Diff}}}
\def\SDiffL{\Sdiff \hskip -.2cm \phantom{a}_{\mathbb{R}}}

\def\proof{\medskip \noindent \textit{Proof: }}

\def\diff{\mathrm{Diff} }
\def\Sdiff{\underset{^{\sim \! \! \sim \! \! \sim}}{\mathrm{Diff}}}

\def\sideremark#1{\ifvmode\leavevmode\fi\vadjust{\vbox to0pt{\vss 
      \hbox to 0pt{\hskip\hsize\hskip1em           
 \vbox{\hsize2cm\tiny\raggedright\pretolerance10000
 \noindent #1\hfill}\hss}\vbox to8pt{\vfil}\vss}}}%
\begin{document}

\title{Peetre-Slov\'{a}k's theorem revisited}


\author{J. Navarro \thanks{Corresponding author. Email address: {\it navarrogarmendia@unex.es} \newline 
Department of Mathematics, Universidad de Extremadura, Avda. Elvas s/n, 06071, Badajoz, Spain. \newline The first author has been partially supported by Junta de Extremadura and FEDER funds. }  
 \and  J. B. Sancho }

\maketitle

\begin{abstract}
In 1960, J. Peetre proved the finiteness of the order of  linear local operators. Later on, J. Slov\'{a}k vastly generalized this theorem, proving the finiteness of the order of a broad class of (non-linear) local operators. 

In this paper, we use the language of sheaves and ringed spaces to prove a simpler version of Slov\'{a}k's result. The statement we prove, adapting Slov\'{a}k's original ideas, deals with local operators defined between the sheaves of smooth sections of fibre bundles, and thus covers many of the applications of Slov\'{a}k's theorem.

\bigskip





\end{abstract}

\tableofcontents

\bigskip

\section*{Introduction}

In 1960, J. Peetre  proved a celebrated theorem on the finiteness of the order of linear local operators; let us briefly recall the statement of this result. 

Let $X$ be a smooth manifold and let $E \to X$, $\bar{E} \to X$ be vector bundles over it. 
Let $ \SDiffL ^k (E , \bar{E})$ denote the sheaf of $\mathbb{R}$-linear differential operators of order $\leq k$ between them, 
and consider the direct limit $$  \SDiffL (E , \bar{E}) := \lim \limits _\rightarrow  \SDiffL ^k (E , \bar{E}) \ , $$
so that elements in the vector space of global sections $\diff _\mathbb{R} (E, \bar{E})$ are linear differential operators locally of finite order.
 

On the other hand, let $\mathcal{E}, \bar{\mathcal{E}}$ be the sheaves of smooth sections of $E, \bar{E}$. A morphism of sheaves $\phi \colon \mathcal{E} \to \bar{\mathcal{E}}$ is $\mathbb{R}$-linear if $\phi (\lambda s + \mu s') = \lambda \phi(s) + \mu \phi(s')$, for any real numbers $\lambda , \mu$ and any sections $s, s'$ defined on any open set. Let $\mathrm{Hom}_{\mathbb{R}} ( \mathcal{E} , \bar{\mathcal{E}} )$ denote the vector space of linear morphisms of sheaves.

Clearly, any linear differential operator $P \colon E \rightsquigarrow \bar{E}$ defines a linear  morphism of sheaves $ \phi_P \colon \mathcal{E} \to \bar{\mathcal{E}}$, and Peetre's theorem affirms the reciprocal: any linear morphism of sheaves is produced by a differential operator.

\begin{Teo}[Peetre, \cite{Peetre}]
The map $P \mapsto \phi_P$ establishes a linear isomorphism:
$$ \diff _{\mathbb{R}} ( E , \bar{E} ) = \mathrm{Hom}_{\mathbb{R}} ( \mathcal{E} , \bar{\mathcal{E}} ) \ . $$
\end{Teo}

Later on, Cahen, De Wilde and Gutt (\cite{Cahen}) proved a multilinear version of this result in the course of their investigations on Hochschild cohomology of smooth manifolds. But it was J. Slov\'{a}k whom, in the context of the theory of natural bundles,  vastly generalized Peetre's theorem, proving the finiteness of the order of a broad class of (non-linear) local operators (\cite{Slovak}). Since then, other authors (\cite{Chrastina}, \cite{Zajtz}) have found various classes of local operators for which different finiteness properties apply.
Nevertheless, all these non-linear Peetre theorems remain rather technical, whereas, as Slov\'{a}k himself pointed out (\cite{SlovakJDG}, Section 1.3), a weaker statement was enough for many of the applications.

In this paper, following Slov\'{a}k's ideas, we give a complete proof of a simple version of this Peetre-Slov\'{a}k's theorem, whose level of generality is enough for the main applications in the theory of natural operations (see recent work in \cite{Einstein}, \cite{Branas}, \cite{Differential}, and compare to those in \cite{Kolar}, \cite{SlovakJDG}), as well as for those in the geometric theory of variational calculus (\cite{KolarVariational}).

To be precise, let $F \to X$, $ \bar{F} \to X$ be fibre bundles over a smooth manifold $X$. 
A (possibly non-linear) differential operator between them is a ``smooth'' map $P \colon JF \to \bar{F}$ over $X$. Here, $JF$ denotes the space of $\infty$-jets of $F$,  endowed with the smooth structure that inherits as the inverse limit, in the category of ringed spaces, of the sequence of $k$-jet prolongations (see details in Section \ref{RingedSection}). 

On the other hand, let $\mathcal{F}, \bar{\mathcal{F}}$ be the sheaves of smooth sections of $F$ and $\bar{F}$. A morphism of sheaves $\phi \colon \mathcal{F} \to \bar{\mathcal{F}}$ is regular if it maps smooth families of sections of $\mathcal{F}$ into smooth families of sections of $\bar{\mathcal{F}}$. 

As an example, any differential operator $P \colon JF \to \bar{F} $ defines a regular morphism of sheaves:
$$ \phi_P \colon \mathcal{F} \to \bar{\mathcal{F}} \quad , \quad s \ \mapsto \ P(j^\infty_x s) \ . $$

This paper is devoted to the proof of the following statement (that is a particular case of \cite{Slovak}, Thm. 2):

\begin{TeoNumerado}
The map $P \mapsto \phi_P$ defines a bijection:
$$ \diff ( F , \bar{F} )   \, = \,  \reg ( \mathcal{F} , \bar{\mathcal{F}} ) \ ,  $$ 
where $\diff ( F , \bar{F} )$ stands for the set of differential operators, and $ \reg ( \mathcal{F} , \bar{\mathcal{F}} )$ for that of regular morphisms of sheaves.
\end{TeoNumerado}

This note is organized as follows: in the first Section, we review the smooth structure of the space of jets, as well as the statement of Whitney's Extension Theorem, that will be used in the subsequent proof. 

The second Section is devoted to the proof of some technical results, of independent interest, and the key Lemma \ref{LemaNonLinearPeetre}. 

Finally, we accomplish the proofs of Theorem \ref{PeetreThm} and Theorem \ref{NonLinearPeetre} in sections three and four, respectively.

\section{Preliminaries}


\subsection{The ringed space of jets}\label{RingedSection}

Let $F \to X$ be a fibre bundle. Its jet prolongations produce a sequence of smooth maps:
$$ \ldots \to J^k F \to J^{k-1}F \to \ldots \to F \to X \ . $$

The $\infty$-jet prolongation $JF$ is intuitively defined as the inverse limit of this sequence. As the category of smooth manifolds does not posses inverse limits, we will work in the larger category of ringed spaces:

\begin{definition}
 A ringed space is a pair $(X , \mathcal{O}_X)$, where $X$ is a topological space and $\mathcal{O}_X $ is a 
sub-algebra of the sheaf of real-valued continuous functions on $X$. 

A morphism of ringed spaces $\varphi \colon (X , \mathcal{O}_X) \to (Y , \mathcal{O}_Y)$ is a continuous map $\varphi \colon X \to Y $ such that composition with $\varphi$ induces a morphism of sheaves:
$$ \varphi^* \colon \mathcal{O}_Y \to \varphi_* \mathcal{O}_X \ . $$ 

The set of morphisms of ringed spaces will be denoted $\mathrm{Hom}(X,Y)$.
\end{definition}

As an example, any smooth manifold $X$ is a ringed space, where $\O_X = \mathcal{C}^\infty_X$ is the sheaf of smooth real-valued functions. If $X$ and $Y$ are smooth manifolds, a morphism of ringed spaces $X \to Y$ is just a smooth map.


\begin{definition}
The space of jets, $JF$, is the ringed space $(JF, \O _J)$ defined as follows:
\begin{itemize}
\item the underlying topological space is the inverse limit of the topological spaces $J^kF$; i.e., is the set: 
$$ JF := \lim_\leftarrow J^kF  $$ 
endowed with the minimum topology for which the natural projections $\pi_k \colon JF \to J^kF $ are continuous.

\item the sheaf of smooth functions is:  
\begin{equation}\label{HazFunciones} \O_J := \lim_\rightarrow \pi^{-1}_k \O_{J^k} \end{equation} 
so that, locally, every smooth function on $JF$ can be written as $h \circ \pi_k$, for some smooth function $h$ on some 
$X_k$.
\end{itemize} 
\end{definition}

With this sheaf $\O_J$ of smooth functions, 
the canonical projections $\pi_k \colon (JF, \O_J) \to (J^kF , \mathcal{O}_{J^k})$ are 
morphisms of ringed spaces.






\begin{theorem}[Universal Property]
For any ringed space $(Y , \O_Y)$, the projections $\pi_k$ induce a bijection, functorially on $Y$:
$$ \mathrm{Hom} (Y, JF) \ = \ \lim_\leftarrow \mathrm{Hom} (Y, J^kF ) \ . $$ 
\end{theorem}

\proof One inclusion being obvious, let us only check that if $\varphi \colon Y \to JF$ is a continuous map such that $\pi_k \circ \varphi$ is a morphism of ringed spaces for any $k \in \mathbb{N}$, then $\varphi$ is a morphism of ringed spaces.

Let $f \in \mathcal{O}_J (U)$ be a smooth function and let $y \in \varphi^{-1}( U)$. On a neighbourhood of $\varphi (y)$, the function $f$ factors as $f = f_k \circ \pi_k$ for some
$f_k \in \mathcal{O}_{J^kF} (V_k)$, and therefore:
$$ \varphi^* f = \varphi^* ( f_k \circ \pi_k) = (\pi_k \circ \varphi)^* f_k $$ that belongs to $ \mathcal{O}_Y ( (\pi_k \circ \varphi)^{-1} V_k) $    because $\pi_k \circ \varphi$ is a morphism of ringed spaces.

\hfill $\square$

\begin{definition}
Let $s \colon X \to F$ be a smooth section. Due to the Universal Property, 
the $k$-jet prolongations $j^k s \colon X \to J^kF$ produce a morphism of ringed spaces, $$ j^\infty s \colon X \to J F \quad , \quad j^\infty s \, (x) = j^\infty_x s \ , $$ called the jet prolongation of $s$.
\end{definition}	


\begin{proposition}\label{MorfismosAVariedad}
Let $Z$ be a smooth manifold. A continuous map $\varphi \colon JF \to Z$ is a morphism of ringed spaces if and only if
for any point $x \in JF $, there exist $k \in \mathbb{N}$, an open set $U \subset J^kF$ and a smooth map $\varphi_k \colon U \to Z$ such that:
$$ 
\xymatrix@C3mm{
\pi_k^{-1} (U) \ar[rr]^-{\varphi} \ar[dr]_{\pi_k} & &   Z  \\
& U \ar[ur]_{\varphi_k} & 
} \ .$$
\end{proposition}

\proof Let $\varphi \colon JF \to Z$ be a continuous map, let $x \in JF $ be a point and let $(V, z_1, \ldots , z_n)$ be a 
coordinate chart around $\varphi (z) $ in $Z$.

If $\varphi$ is a morphism of ringed spaces, each of the  functions $z_1 \circ \varphi, \ldots , z_n \circ \varphi \in \mathcal{O}_\infty (\varphi^{-1} V)$ locally factors  through
some $J^lF$. As they are a finite number, there exists $k \in \mathbb{N}$ and an open neighbourhood $U$ of $x$ such that all of them, when restricted to $U$,
factor through $J^kF$. Hence, $\varphi_{|U} = (\varphi_k \circ \pi_k)_{|U}$, where $\varphi_k \equiv ( z_1 \circ \varphi , \ldots , z_n \circ \varphi)$.

Conversely, if $\varphi = \varphi_k \circ \pi_k$ on a neighbourhood of a point, then $\varphi^* f = f \circ \varphi = f \circ ( \varphi_k \circ \pi_k) = 
((f \circ \varphi_k) \circ \pi_k)$ on that neighbourhood, so that $\varphi^* f$ is a smooth function on $JF$ and $\varphi$ is a morphism of ringed spaces.
\hfill $\square$

\subsection{Whitney's Extension Theorem}

For any pair of multi-indexes, 
$$ I = (r_1, \ldots , r_n) \ , \  J = (\overline{r}_1 , \ldots , \overline{r}_n) \ \in \, \mathbb{Z}^+ \times \stackrel{n}{\ldots} \times \mathbb{Z}^+ \ $$ 
let us introduce the following notations:
$$ | I | := r_1 + \ldots + r_n \quad , \quad I! := r_1 ! \ldots r_n! \quad , \quad  I + J := ( r_1 + \overline{r}_1 , \ldots , r_n + \overline{r}_n) \ . $$ 

Some other standard multi-index notation will be used throughout this Section, without more explicit mention. As an example, the Taylor expansion of a smooth function $f$ on $\mathbb{R}^n$ at a point $a \in \mathbb{R}^n$ is denoted:
$$ \mathtt{T}_a f := \sum_{I } \, \frac{D_I f (a) }{I!} \ (x-a)^I \ . $$ 

\begin{definition}
A Taylor expansion $\mathtt{T}_a$ at a point $a\in \mathbb{R}^n$ is an arbitrary series:
$$ \mathtt{T}_a = \sum_{I } \, \frac{\lambda_{I,a}}{I!} \ (x-a)^I \qquad \lambda_{I,a} \in \mathbb{R} \ . $$
A family of Taylor expansions $\{ \mathtt{T}_a \}_{a\in K}$ on a set $K \subset \mathbb{R}^n$ defines functions: 
$$ \lambda_I \, \colon \, K \to \mathbb{R} \quad , \quad  a \, \mapsto \, \lambda_{I,a} \ . $$  
\end{definition}

Let $K \subset \mathbb{R}^n$ be a compact set and consider a family of Taylor expansions $\{ \mathtt{T}_a \}_{a\in K}$ on the points of $K$.
This Section deals with the question of whether there exists a smooth function $f \in \mathcal{C}^\infty (\mathbb{R}^n)$ such that $\mathtt{T}_a f = \mathtt{T}_a $, for all $a \in K$.

A necessary condition is given by Taylor's Theorem: if $f$ is a smooth function on $\mathbb{R}^n$ and $K \subset \mathbb{R}^n$ is a compact set,  then for any  $m \in \mathbb{N}$ and $\epsilon >0$, there exists $\delta > 0 $ such that:
\begin{equation}\label{TaylorTheorem}
 x,y \in K \ , \  \| x - y \| < \delta \quad \Rightarrow \quad \left| \, f(y) - \sum_{|J|=0}^m \frac{D_J f(x)}{J!} (y-x)^J \, \right| \leq \epsilon \, \| y-x \|^m \ .
\end{equation}

Whitney's Theorem provides a sufficient condition:  

\medskip
\noindent {\bf Whitney's Extension Theorem (\cite{Whitney}):} {\it 
 Let $K \subset \mathbb{R}^n$ be a compact set and let $\{ \mathtt{T}_a \}_{a \in K}$ be a family of Taylor expansions on $K$.
$$ \mbox{ There exists } f \in \mathcal{C}^\infty (\mathbb{R}^n) \mbox{ such that } \ \mathtt{T}_a f = \mathtt{T}_a \mbox{ for any }  a \in K  $$ if and only if the following condition (that we will refer to as {\rm "Taylor's condition"}) holds:
$$ \mbox{For any given I, $m \in \mathbb{N}$, and $\epsilon > 0$, there exists $\delta > 0$ such that:} $$
\begin{equation*}\label{TaylorCondition}
x,y \in K \ , \ \| x - y \| < \delta \quad \Rightarrow \quad \left| \, \lambda_{I ,y} - \sum_{|J| =0 }^m  \frac{\lambda_{I+J,x}}{J!} \ (y-x) ^J \, \right| \, \leq \, \epsilon \, \|y-x \|^m \ .
\end{equation*}  }

\section{Morphisms of sheaves} 

This Section is devoted to the proof of Propositions {\ref{MorfismoArbitrario} and \ref{LemaDeLema}, regarding arbitrary morphisms of sheaves between sheaves of smooth sections of fibre bundles.

\medskip

\begin{lemma}\label{PropoCono}
Let $K \subset \mathbb{R}^n$ be the truncated 
cone on $\mathbb{R}^n$ defined by the equations:
\begin{equation}\label{Cone}
  K:= \{ x \in \mathbb{R}^n \colon \, x_1^2 + \ldots + x_{n-1}^2 \leq x_n^2 \quad , \quad | x_n | \leq 1 \ \} \  \ 
\end{equation} 
and let $K_1 := K \cap \{ x_n \geq 0 \}$ and $ K_2 := K \cap \{ x_n \leq 0 \}$.

If $u , v \in \mathcal{C}^\infty (\mathbb{R}^n)$ are smooth functions with the same Taylor expansion at the origin, $\mathtt{T}_0 u = \mathtt{T}_0 v $, then there exists $ f \in \mathcal{C}^\infty   (\mathbb{R}^n)$ such that $$  f_{|K_1} = u_{|K_1} \quad  , \quad f_{|K_2} = v_{|K_2}  \ .  $$
\end{lemma}

\proof It is enough to argue the case $v = 0$. To do so, let us check we can apply Whitney's  Theorem to the following  family of Taylor expansions:
$$ \mathtt{T}_a  :=  \begin{cases} \mathtt{T}_a u \ , \ \mbox{ if } a \in K_1 \\ 
\quad  0 \ \,  , \, \ \mbox{ if } a \in K_2
\end{cases} \ . $$

Taylor's condition trivially holds
whenever $x,y \in K_1$ or $x,y \in K_2$.
If $x \in K_1$ and $y \in K_2$ (the case where $x$ and $y$ are interchanged is analogous), then: 
$$ \| x \| , \| y \| \leq \| y-x \| \ . $$ 

Given $I, m , \epsilon $, the hypothesis $\mathtt{T}_0 u = 0 $ implies  there exists $\delta > 0$ such that, for any $J$ with $|J| \leq m$,
$$ x \in K_1 \ , \ \| x \| \leq \delta \quad \Rightarrow \quad | \lambda_{I+J , x} | \leq \, \epsilon \, \| x \|^m \ . $$ 

A smaller $\delta$, if necessary, also guarantees $| (y-x)^J | < 1 $ whenever $\| x - y \| < \delta$ and $|J|\leq m$.

Therefore, whenever $x \in K_1 $ and $y \in K_2 $ satisfy $\| x - y \| \leq \delta$,
\begin{align*}
 \left| \, \lambda_{I,y}  - \sum_{|J| \leq m} \frac{\lambda_{I+J , x} }{J!} (y-x)^J \, \right| \, & \leq \, 0  +   \sum_{|J| \leq m}  \left| \, \lambda_{I+J , x} \, \right|  \, \frac{| ( y-x )^J |}{J!} \, \\
 &  \leq  \ \epsilon \, \| x \|^m \left(  \sum_{|J|=0}^m \frac{1}{J!} \right) \, \leq \ \overline{\epsilon}   \| y-x \|^m   .
 \end{align*} 
\hfill $\square$

\begin{proposition}\label{MorfismoArbitrario} Let $\phi \colon \mathcal{F} \to \overline{\mathcal{F}}$ be a  morphism of sheaves. 
For any sections $s,s'$ of $\mathcal{F}$ defined on a neighbourhood of a point $x \in X$:
$$ j^\infty_x s = j^\infty_x s' \qquad \Rightarrow \qquad \phi(s)(x) = \phi(s')(x) \ . $$
\end{proposition}  

\medskip
\noindent \textit{Proof:} As the statement  is local, we can suppose $x=0$ is the origin of $X=\mathbb{R}^n$ and $F = \mathbb{R}^r \times \mathbb{R}^n$ is trivial. 
We can also assume $\overline{F} = \mathbb{R} \times \mathbb{R}^n$ is trivial, with one-dimensional fibres.

%

Hence, let $s \equiv (s_1 , \ldots , s_r), s' \equiv (s'_1, \ldots , s'_r) \colon \mathbb{R}^n \to \mathbb{R}^r$ be smooth maps with the same $\infty$-jet at the origin.
Let $K = K_1 \cup K_2 \subset \mathbb{R}^n$ be a truncated cone as in 
(\ref{Cone}).

 In this situation, Lemma \ref{PropoCono} proves the existence of smooth functions $f_1, \ldots f_r $ on $\mathbb{R}^n$
such that:
$$ f_{i|K_1} = s_{i|K_1} \quad , \quad f_{i|K_2} = s'_{i|K_2} \quad  \quad i = 1, \ldots m \  . $$

The section $f \equiv (f_1 , \ldots , f_r)$ satisfies $ f_{|K_1} = s_{|K_1}$, $f_{|K_2} = s'_{|K_2}$, and, consequently, 
as $\phi$ commutes with restrictions to open sets:
$$ \phi (f)_{|\stackrel{\circ}{K_1}} = \phi (s)_{|\stackrel{\circ}{K_1}} \qquad , \qquad \phi(f)_{|\stackrel{\circ}{K_2}} = \phi(s')_{|\stackrel{\circ}{K_2}} \ . $$ 

By continuity,
$$ \phi(s) (x) \, = \, \phi(f) (x) \, = \, \phi(s') (x) \ . $$ \qed

\medskip

As a consequence, if  $\phi \colon \mathcal{F} \to \overline{\mathcal{F}}$ is a  morphism of sheaves,  the following map (between sets) is well-defined:
$$P_{\phi} \colon J^\infty F \to \overline{F} \quad , \quad P_{\phi} (j^\infty_x s) := \phi (s) (x) \ . $$ 

\begin{proposition}\label{LemaDeLema}
Let $\phi \colon \mathcal{F} \to \overline{\mathcal{F}}$ be a  morphism of sheaves.
If $s, s'$ are sections of $\mathcal{F}$ and $x_k \to x$ is a sequence converging to a point $x \in X$, then, 
$$ j^k_{x_k} s = j^k_{x_k} s' \ , \ \forall \ k \in \mathbb{N} \quad \Rightarrow \quad  \exists \ k_0 \ \colon \ 
\phi(s)(x_k ) = \phi(s') (x_k) \quad \forall \ k > k_0 \ . $$
\end{proposition}

\proof Again, we can suppose $X = \mathbb{R}^n$ and $F= \mathbb{R}^r \times \mathbb{R}^n$, $\overline{F} = \mathbb{R} \times \mathbb{R}^n$ are trivial bundles. 

If the statement is not true, taking a subsequence we can assume there exists $s, s'$ such that, for any $k$:
$$ j^k_{x_k} s = j^k_{x_k} s' \quad , \quad \phi(s) (x_k) \neq \phi(s') (x_k) \ . $$




Let us consider another sequence $y_k \to x$ such that:
\begin{align}\label{CondicionSmooth} 
| \phi(s)(x_k) - \phi(s')(y_k) | \, &> \, k \, \| y_k - x_k \| \ . \\ 
x_k &\neq y_l \qquad \forall \, k,l \ . 
\end{align}


Let us now apply Whitney's Extension Theorem on the compact $(x_k) \cup (y_k) \cup \{ 0 \}$ to the family of jets:
\begin{equation}\label{FamilyOfJets}
 \{ \  j^\infty_{x_k} s \ , \ j^\infty_{y_k} s' \ , \ j^\infty_x s = j^\infty_x s' \ \} \ . 
\end{equation} 

Due to Taylor's Theorem (\ref{TaylorTheorem}), given $I,m, \epsilon$, there exists $\delta > 0$ such that, for any $i=1, \ldots , m$ and any points $a, b$ in the compact, the condition $\| a - b \| < \delta $ implies:
\begin{align}
&\left| (D_I s_i)(a) - \sum_{|J|=0}^m (D_{I+J} s_i)(b) \frac{(a-b)^J}{J!} \right| \ \leq \ \epsilon \, \| a - b \|^m \\
\label{TaylorU}
&\left| (D_I s'_i)(a) - \sum_{|J|=0}^m (D_{I+J} s'_i)(b) \frac{(a-b)^J}{J!} \right| \ \leq \ \epsilon \, \| a - b \|^m 
 \end{align} 


A smaller $\delta$, if necessary, also guarantees:
$$ \| x_k - y_l \| < \delta \quad \Rightarrow \quad k,l \, > \, |I| + m =: M \ , $$
and this allows  to invoke Whitney's Theorem:

\begin{itemize}
\item[-] If $\| x_k - x_l \| , \|  y_k - y_l \| < \delta$, then (\ref{TaylorU}) is just Taylor's condition.

\item[-] If $ \| x_k - y_l \| < \delta$, then $k , l > M = |I| + m $; as $j^M_{x_k} s = j^M_{x_k} s'$, it holds:
\begin{align*}
 \left| (D_I s') (y_l) -  \sum_{|J|=0}^m  (D_{I+J}s)(x_k) \frac{(y_k-x_k)^J}{J!} \right| = \left| (D_I s') (y_l) - \sum_{|J|=0}^m (D_{I+J} s')(x_k) \frac{(y_k-x_k)^J}{J!} \right| 
\end{align*}
and inequality (\ref{TaylorU}) proves this quantity is less or equal than $\epsilon \, \| y_l - x_k \|^m $.
\end{itemize}

Therefore, there exists a smooth map $f \colon \mathbb{R}^n \to \mathbb{R}^r$ realizing the family of jets (\ref{FamilyOfJets}). 
Due to Proposition \ref{MorfismoArbitrario}, this map $f$ satisfies:
$$ \phi(f)(x_k) = \phi(s) (x_k) \quad , \quad \phi(f)(y_k) = \phi(s')(y_k) $$ and (\ref{CondicionSmooth}) contradicts the smoothness of $\phi(f)$. 

\qed

\subsection{Linear morphisms of sheaves}

In this Section we prove Proposition \ref{LemmaPetreeLineal}, which deals with linear morphisms of sheaves, and is the 
main ingredient in Peetre's theorem.

\begin{lemma}\label{LemaSucesion}
 Let $a_k \to 0$ be a sequence of 
points converging to the origin of $\mathbb{R}^n$. 

Assume there exist $c\in \mathbb{N}$ such that:
$$\| a_k \| , \| a_l \| \, < \, c \| a_k - a_l \| \qquad \forall \, k \neq l \ , $$ 
and let $\{ \mathtt{T}_{a_k} , \mathtt{T}_0 = 0  \}_{k \in \mathbb{N}}$ be a 
family of Taylor expansions on $\, \{ a_k \}, 0$. 

If Taylor's condition holds at the origin; i.e., for any $I$, $m$, $\epsilon$, there exists
$\nu \in \mathbb{N}$ such that :
$$ k > \nu \quad \Rightarrow \quad |\lambda_{I, a_k} | \, \leq \, \epsilon \, \| a_k \| ^m \  $$
then there exists a global smooth function $f \in \mathcal{C}^\infty (\mathbb{R}^n)$ realizing those expansions:
$$ \mathtt{T}_{a_k} f  = \mathtt{T}_{a_k}  \qquad , \qquad \mathtt{T}_{0} f = 0 \ . $$ 
\end{lemma}

\proof We can assume $| (a_k - a_l)^J | < 1$, for any multi-index $J$ and any $k,l$.


Given $I, m, \epsilon $, the hypothesis allows to find $\nu \in \mathbb{N}$ such that, for any $J$ with $|J| \leq  m $:
$$ k > \nu \quad \Rightarrow \quad |\lambda_{I+J, a_k} | \, \leq \, \epsilon \, \| a_k \|^m \ .   $$ 

Hence, for any $k,l > \nu$:
\begin{align*}
 \left| \, \lambda_{I,a_k} -   \sum_{|J|=0}^m \lambda_{I+J , a_l}\, \frac{(a_k - a_l)^J}{J!} \, \right| &\leq   \, |\lambda_{I , a_k} | + \sum_{|J|=0}^m |\lambda_{I+J, a_l}| \, \frac{|(a_k - a_l)^J|}{J!} \\
& \leq \, \epsilon \, \| a_k \|^m + \, \epsilon  \, \| a_l \|^m \sum_{|J|=0}^m \frac{ 1 }{J!} \\
& \leq
\, \epsilon \left( c^m \| a_k - a_l \|^m + c^m M \| a_k - a_l \|^m \right)  \leq \,   \overline{\epsilon} \, \| a_k - a_l \|^m  \ . 
\end{align*}

Therefore, Whitney's Theorem applies and the thesis follows.

\qed

\bigskip

Let $\,E , \bar{E} \to X\,$ be vector bundles over a smooth manifold $\,X$, and let $\,\mathcal{E}, \bar{\mathcal{E}} \, $ denote their sheaves of smooth sections, respectively.

\begin{proposition}\label{LemmaPetreeLineal} Let $\phi \colon \mathcal{E} \to \bar{\mathcal{E}}$ be an $\mathbb{R}$-linear morphism of sheaves.

For any point $x \in X$, there exist an open
neighborhood $U \subset X$ and a natural number $k \in \mathbb{N}$ such that:
$$ y \in U \quad , \quad j^k_y s_1 = j^k_y s_2 \qquad \Rightarrow \qquad \phi(s_1)(y) = \phi(s_2) (y) \ , $$ where $s_1$ and $s_2$ are any representatives of $\,j^k_y s_1\,$ and $\,j^k_y s_2$.
\end{proposition}

\medskip \noindent \textit{Proof:} We can assume $\,X= \mathbb{R}^n$, $\,\mathcal{E} = \mathcal{O}_{\mathbb{R}^n}\,$ and $x $ is the origin. By linearity, it is enough to prove that  there exist an open neighborhood $U \subset \mathbb{R}^n$ of the origin
and $\,k \in \mathbb{N}\,$ such that:
\begin{equation}\label{ToProve} y \in U \quad , \quad j^k_y f = 0 \qquad \Rightarrow  \qquad \phi (f) (y) = 0 \ .
\end{equation}

\medskip
To do so, we first prove a weaker statement; namely, that there exist a constant $\,M > 0$, an open neighborhood $\,U\,$ of the origin
and $\,k \in \mathbb{N}\,$ such that:
\begin{equation}\label{Simplified}   y \in U - \{ 0\} \quad , \quad j^k_y f = 0 \qquad \Rightarrow  \qquad |\phi (f) (y)| < M \ . \end{equation}

If this were not true, considering the neighbourhoods $\,U_m := \{ \| x \| \leq \frac{1}{m} \} $, we can produce a sequence of different points $x_m \to 0$, $x_m \in U_m - \{ 0 \}$ and functions $\,h_m\,$ such that:
$$ j^k_{x_m} h_m = 0 \quad \mbox{ but } \quad |\phi(h_m) (x_m) | > m \ . $$ 

By Lemma \ref{LemaSucesion}, there exists a global smooth function $f$ such that:
$$ j^\infty_{x_k} f = j^\infty_{x_k} h_k \qquad , \qquad j^\infty_0 f = 0 \ . $$ 

These conditions imply $\phi(f)(x_k) = \phi(h_k)(x_k)$ and $\,\phi(f)(0)=0$ (see Lemma \ref{MorfismoArbitrario}), so we arrive to contradiction: 
$$ \phi(f) (x_k ) \to \phi(f) (0) = 0 \qquad  \mbox{ while } \qquad | \phi (f) (x_k) | = | \phi (h_k) (x_k) | > k \ . $$

\medskip
Now, observe that (\ref{Simplified}) implies (\ref{ToProve}) on any point but the origin, for if there exist $\,y \in U - \{ 0 \}\,$ and a a function $\,f\,$ such that $\,j^k_y f = 0\, $ but $\,|\phi(f)(y)| = \epsilon >0$, then, rescaling the function $\,f\,$ by a factor $\,2 M / \epsilon$, 
we have:
$$ j^k_y \left( \frac{2M}{\epsilon} \, f \right) = 0 \quad \mbox{ but } 
\quad \left| \phi\left( \frac{2 M}{\epsilon}\, f \right) (y) \right| = \frac{2M}{\epsilon}\, | \phi(f)(y)| = 2 M > M \ .$$

\medskip
Finally, to prove (\ref{ToProve}) on the origin we argue as follows: let $\,f\,$ be a function such that $\,j^k_0 f = 0$. Consider any sequence of different points $\,x_m \to 0\,$  such that $\,2\, \| x_m - x_{l} \| > \| x_m \| , \| x_l \|$, for $\,m< l$. Let us also consider the sequence of jets $\,\mathtt{T}_{x_m} = \sum_{I } \, \frac{\lambda_{I,x_m}}{I!} \ (x-x_m)^I \,$ defined by the following conditions:
$$ \lambda_{I,x_m} := \left\{ 
\begin{array}{ccl}
(D_I f)(x_m) & , & \mbox{ if } |I| \leq k. \\
0 &, & \mbox{ if } |I| > k.
\end{array} 
\right. $$

Again, we can apply Lemma \ref{LemaSucesion}, that assures the existence of a global smooth function $\,u \,$ such that:
$$ j^\infty_{x_m} u = \mathtt{T}_{x_m} \quad , \quad j^\infty_0 u =0\ . $$

Due to (\ref{Simplified}), the condition $j^k_{x_m} u = j^k_{x_m} f$ implies $\phi(u) (x_m) = \phi (f) (x_m) $. On the other hand, $\,\phi (u) (0) = 0$ (Lemma  \ref{MorfismoArbitrario}), so, by continuity, 
$$\phi (f) (0) = \phi(u)(0) = 0 \ . $$

\hfill $\square$

\subsection{Regular morphisms of sheaves}

Let $\, F \to X\,$ be a fibre bundle. 

Given a smooth manifold $\,T$, let us denote $\,X_T := T \times X$. Any open set $\,U \subset X_T\,$ can be thought as a family of open sets $\,U_t \subset X$, where $\,U_t\,$ is the fibre of $\,U \to T\,$ over $\,t \in T$. 

A family of sections $\,\{ \, s_t \colon U_t \to F \, \}_{t \in T}\,$ defines a map: 
$$ s \colon U \to F \quad , \quad s (t,x) := s_t (x) \  $$ and $\,\{ s_t \}_{t \in T}\,$ is said smooth (with respect to the parameters $\,t \in T$) precisely when the map $\,s \colon U \to F\,$ is  smooth:

\begin{definition}
A smooth family of sections of $\,F\,$ parametrized by $\,T\,$ is a section of $\,F\,$ with support on an open set $\,U\,$ of $\,X_T$:
$$
\xymatrix@C=3mm{
 & & F \ar[d] \\
X_T \supset U \ \ar[urr]^{s} \ar[rr] &  & X  }   $$
\end{definition}

\begin{example}\label{UniversalFamily} Let $X = \mathbb{R}^n$, $F = \mathbb{R}^r \times \mathbb{R}^n $ and 
let $T := Pol_k( \mathbb{R}^n , \mathbb{R}^r) $ be the smooth manifold of polynomial maps $f \colon \mathbb{R}^n \to \mathbb{R}^r$ of degree less or equal than $k$. 

The universal family $\xi$ is defined in  $U = J^k F = Pol_k(\mathbb{R}^n , \mathbb{R}^r) \times \mathbb{R}^n$ by the formula:
$$ \xi \colon U \to F = \mathbb{R}^r \times \mathbb{R}^n \quad , \quad (f , x) \ \mapsto \ (f(x) , x) \ ;$$ 
that is, $\xi_f \colon \mathbb{R}^n = X \to F = \mathbb{R}^r \times \mathbb{R}^n$ is the section defined by the polynomial $f$.
\end{example}

\begin{definition}\label{DefinitionRegular}
 A morphism of sheaves $\phi \colon \mathcal{F} \to \overline{\mathcal{F}} $ is regular if, for any smooth family of sections $\{ s_t \colon U_t \to F \}_{t \in T}$, the family $\{ \phi(s_t) \colon U_t \to \overline{F} \}_{t\in T}$ is also smooth. 
\end{definition}




Let $\phi \colon \mathcal{F} \to \mathcal{\bar{F}}$ be a regular morphism of sheaves. As on any morphism of sheaves (Proposition \ref{MorfismoArbitrario}), the following map is well-defined:
$$P_{\phi} \colon J F \to \bar{F} \quad , \quad P_{\phi} (j^\infty_x s) := \phi (s) (x) \ . $$

\begin{lemma}\label{LemaNonLinearPeetre}
If $\phi \colon \mathcal{F} \to \overline{\mathcal{F}}$ is a regular morphism of sheaves, then the map $P_{\phi} \colon J^\infty F \to \overline{F}$ locally factors through some finite jet space. 

That is, for any $j^\infty_x s \in J^\infty F$ there exist an open neighbourhood $V $, a natural number $k \in \mathbb{N}$ and a commutative triangle of maps:
$$ 
\xymatrix{
V \ar[rr]^-{P_{\phi}} \ar[dr]_{\pi_k} & & \overline{F}  \\
& J^k F \ar[ur] & 
} \ .
$$
\end{lemma}

\proof
It is a local problem, so we can assume $X= \mathbb{R}^n$, $F = \mathbb{R}^r \times \mathbb{R}^n$ and $x=0$. Moreover, we can assume  the section $s \colon \mathbb{R}^n \to \mathbb{R}^r$ representing the jet $j^\infty_0s$ is the zero section $s = 0$. 

For each $k \in \mathbb{N}$, consider the following neighbourhood of $j^\infty_0 s$:
$$ U_k := \left[\,  j^\infty_{x} f \in J F \ / \  \ \  \| x \| \leq \frac{1}{2^k} \ , \ | D_If_i (x) | \leq  \left( \frac{1 }{\, 2^k \,} \right)^k \ , \  \, \forall \, |I| \leq k \  , \ i= 1, \ldots , m \  \right] \ .$$

If the thesis is not true, for any $k \in \mathbb{N}$ we can find $ j^\infty_{x_k} f^k , \,  j^\infty_{x_k} h^k \in U_k$ such that $$ j^k_{x_k} f^k = j^k_{x_k} h^k \quad \mbox{ but } \quad \phi(f^k)(x_k) \neq \phi(h^k)(x_k) \ .$$ 

If we could find smooth sections $f,h$ of $F$ such that $ j^\infty_{x_k} f = j^\infty_{x_k} f^k$, $ j^\infty_{x_k} h = j^\infty_{x_k} h^k$, then the above statement contradicts Proposition \ref{LemaDeLema}, as it would be, for any $k \in \mathbb{N}$: $$ j^k_{x_k} f = j^k_{x_k} h \quad \mbox{ but } \quad \phi(f) (x_k) = \phi(f^k) (x_k) \neq \phi(h^k) (x_k) = \phi(h)(x_k) \ . $$ 

But, in order to apply Whitney's  Theorem and extend $j^\infty_{x_k} f^k , \, j^\infty_{x_k} h^k$, the points $x_k$ may be inconveniently placed. 
To overcome this difficulty, let us parametrize by $\mathbb{R}$
and consider the points $$ z_k := \left( \frac{1}{2^k} \ , \ x_k \right) \in \mathbb{R} \times \mathbb{R}^n \ , $$ that satisfy:
$$
\frac{1}{2^k} \leq \| z_k \| \leq \frac{1}{2^{k-1}} \qquad , \qquad  \| z_k \|  , \| z_l \| \, < \, 4 \| z_k - z_l \| \qquad \forall \, k \neq l .
$$ 

Let us extend each section $s$ of $F$ over $X= \mathbb{R}^n$ to the constant family $s (t,x) := s(x)$ over $\mathbb{R} \times X = \mathbb{R} \times \mathbb{R}^n$. 

For any fixed $I , m $, a sufficiently large $k$ assures:  
$$ |D_I f^k_i (z_k) | \, = \, |D_I f^k_i (x_k) | \, < \, \left( \frac{1}{2^k} \right)^k \leq \, \left( \frac{1}{2^k} \right)^{k-m}  \| z_k \|^m \, = \, \epsilon \, \| z_k \|^m \ $$ and idem for the $h_i^k$.


In this situation, Lemma \ref{LemaSucesion} proves there  exist  smooth sections $f = ( f_1 , \ldots , f_r)$, $h = (h_1, \ldots , h_r)$ such that, for any $k$:
$$ j^\infty_{z_k} f = j^\infty_{z_k} f^k \quad , \quad j^\infty_{z_k} h = j^\infty_{z_k} h^k \quad , \quad \forall \,  k \in \mathbb{N} \ .  $$ 


Thus, these sections $f,h$ of $F$ satisfy:
$$ j^k_{z_k} f = j^k_{z_k} h \quad \mbox{ but } \quad \phi(f)(z_k) \neq \phi(h)(z_k) \quad , \quad \forall \,  k \in \mathbb{N} , $$ in contradiction with Proposition \ref{LemaDeLema}. 

\hfill $\square$

\section{Peetre's characterization of linear differential operators}

Let $\,E , \bar{E} \to X\,$ be vector bundles over a smooth manifold $\,X$. 

The jet spaces $J^kE \to X$ are also vector bundles, and hence the ringed space $JE \to X$ inherits a $\mathbb{R}$-linear structure on fibres.

\begin{definition} A linear differential operator
$\,E \rightsquigarrow \bar{E}\,$ is an $\, \mathbb{R}$-linear morphism of ringed spaces over $X$:
$$ P \colon JE \to \bar{E} \ . $$
\end{definition}

Let $\,\mathcal{E}, \bar{\mathcal{E}}\,$ be the sheaves of smooth sections of $\,E\,$ and $\,\bar{E}$. Any linear differential operator $\,P \colon JE \to \bar{E}\,$ defines an $\,\mathbb{R}$-linear  morphism of sheaves:
$$ \phi_P \colon \mathcal{E} \to \bar{\mathcal{E}} \quad , \quad s \ \mapsto \ P(j^\infty_xs) \ . $$

\begin{theorem}[Peetre]\label{PeetreThm}
The map $P \mapsto \phi_P$ defines a linear isomorphism:
$$ \diff _{\mathbb{R}} ( E , \bar{E} ) = \mathrm{Hom}_{\mathbb{R}} ( \mathcal{E} , \bar{\mathcal{E}} ) \ . $$
\end{theorem}

\medskip
\noindent \textit{Proof:} 
Let $\phi \colon \mathcal{E}_{|V} \to \bar{\mathcal{E}}_{|V}$ be a $\mathbb{R}$-linear morphism of sheaves over an open set $V \subset X$, and let $x \in V$ be a point. We have to check that, on a sufficiently small neighbourhood of $x$, the morphism $\phi$ is defined by a linear differential operator of finite order.

By Proposition \ref{LemmaPetreeLineal}, there exists a neighbourhood $U \subset V$ of $x$ and a natural number $k \in \mathbb{N}$ such that the following map is well-defined:
$$ P \colon (J^k E)_U \to \bar{E}_U \quad , \quad j^k_y s \ \mapsto \ \phi(s) (y) \  $$ where $s$ is any representative of the jet $j^k_y s$.

The map $P$ is $\mathbb{R}$-linear, and we want to prove it is smooth. To do so, we can assume $E$ and $\bar{E}$ are trivial bundles of rank one over $\mathbb{R}^n$, and therefore:
$$ P = \sum_{|I|=0}^k P^I (x_1, \ldots x_n) D_I \ .$$ 

The smoothness of the functions $P^I$ is easily achieved by induction, 
using that $P \circ j^k f = \phi $ is smooth, for any $k$-jet prolongation $j^k f$. 


\hfill $\square$

\section{Slov\'{a}k's characterization of differential operators}

Again, let $F \to X$, $\bar{F} \to X$ denote fibre bundles and let $\mathcal{F}$, $\bar{\mathcal{F}}$ be their sheaves of smooth sections.

\begin{definition} A differential operator 
$\,F \rightsquigarrow \bar{F}\,$ is a morphism of ringed spaces over $X$:
$$ P \colon JF \to \bar{F} \ . $$


\end{definition}

Any differential operator $P \colon JF \to \bar{F}$
allows to define a morphism of sheaves:
$$ \phi_P \colon \mathcal{F} \to \bar{\mathcal{F}} \qquad , \qquad \phi_P (s)(x) := P (j^\infty_x s ) \ .
 $$

This morphism of sheaves is regular, for its value on a smooth family of sections $s_t \colon U_t \to F$ is:
$$ \phi_P (s_t) \ \colon \ (t,x) \ \mapsto \ (t ,  P (j^\infty_x s_t) ) \   $$ that is also smooth on the parameters $t \in T$.


\begin{theorem}\label{NonLinearPeetre}   
Let $\mathcal{F}, \bar{\mathcal{F}}$ be the sheaves of smooth sections of two fibre bundles $F, \bar{F}$ over a smooth manifold $X$.

The map $P \mapsto \phi_P$ establishes a bijection:
$$  \reg ( \mathcal{F} , \bar{\mathcal{F}} ) \, = \, \diff ( F , \bar{F} )  \  $$ where $\reg ( \mathcal{F} , \bar{\mathcal{F}} ) $ and $\diff ( F , \bar{F} )$ stand for the sets of regular morphisms of sheaves and differential operators, respectively.
\end{theorem}

\proof Let $\phi \colon \mathcal{F} \to \overline{\mathcal{F}}$ be a regular morphism of sheaves. 


In virtue of Lemma \ref{LemaNonLinearPeetre},
any point in $J^\infty F$ has an open neighbourhood $V$ such that the following triangle commutes:
$$ 
\xymatrix{
V \ar[rr]^-{P_{\phi}} \ar[dr]_{\pi_k} & & \overline{F}  \\
& J^k F \ar[ur]_{P_k} & 
} \ .$$

It only rests to check that $ P_k \colon J^k F \to \overline{F}$, $j_x^k s \mapsto \phi(s) (x)$ is smooth, and, to this end, we can suppose $X = \mathbb{R}^n$ and $F = \mathbb{R}^r \times \mathbb{R}^n $ is trivial.

Let $\xi$ be the universal family of $k$-jets (Example \ref{UniversalFamily}). 
The smoothness of $P_k$ follows from the regularity of $\phi$, because:
$$ P_k (j^k_xs )  \, = \, \phi (\xi_f)(x)   \ , $$ where $f \colon X \to F$ is the only polynomial whose $k$-jet at $x$ is $j^k_x s$.

\qed


\end{document}